\documentclass[12pt]{amsart}
\usepackage{amsmath, amssymb, latexsym, amsthm, bm}
\usepackage{mathrsfs}
\usepackage[all]{xy}

\newcommand{\Pa}[9]{\bibitem{#1} {#2}, \emph{#3}, {#4} \textbf{#5} ({#6}), {#7}--{#8}.}

\newcommand{\ed}{
\end{document}
}

      \newenvironment{changemargin}[2]{\begin{list}{}{
         \setlength{\topsep}{0pt}\setlength{\leftmargin}{0pt}
         \setlength{\rightmargin}{0pt}
         \setlength{\listparindent}{\parindent}
         \setlength{\itemindent}{\parindent}
         \setlength{\parsep}{0pt plus 1pt}
         \addtolength{\leftmargin}{#1}\addtolength{\rightmargin}{#2}
         }\item }{\end{list}}

\newcommand{\arx}[3]{\texttt{http://arxiv.org/#1/#3}}
\newcommand{\Arh}{Arhangel'ski\u{\i}}
\newcommand{\znb}[2]{{\bm{[}#1; #2\bm{]}}}
\newcommand{\zero}{\mathbf{0}}

\newcommand{\alephes}{{\aleph_0}}

\newcommand{\Dfin}{\mathfrak{D}_\mathrm{fin}}

\newcommand{\bq}{\begin{quote}}
\newcommand{\eq}{\end{quote}}
\newcommand{\cl}[1]{\overline{#1}}

\newcommand{\inv}{^{-1}}

\newcommand{\N}{\mathbb{N}}
\newcommand{\NN}{{\N^{\N}}}

\newcommand{\op}{\operatorname}

\newcommand{\scrA}{\mathscr{A}}
\newcommand{\scrB}{\mathscr{B}}

\newcommand{\cF}{\mathcal{F}}
\newcommand{\cS}{\mathcal{S}}

\newcommand{\cN}{\mathcal{N}}
\newcommand{\cO}{\mathcal{O}}

\newcommand{\cU}{\mathcal{U}}
\newcommand{\cC}{\mathcal{C}}
\newcommand{\Union}{\bigcup}
\newcommand{\cV}{\mathcal{V}}

\long\def\forget#1\forgotten{}

\newcommand{\fd}{\mathfrak{d}}

\newcommand{\oo}{\infty}

\newcommand{\w}{\omega}
\newcommand{\x}{\times}

\newcommand{\nin}{\not\in}

\newcommand{\sbst}{\subseteq}
\newcommand{\spst}{\supseteq}
\newcommand{\sm}{\setminus}
\newcommand{\rest}{\restriction}

\newcommand{\cov}{\op{cov}}

\newtheorem{thm}{Theorem}%[section]
\newcommand{\bthm}{\begin{thm}}
\newcommand{\ethm}{\end{thm}}

\newcommand{\bprp}{\begin{prp}}
\newcommand{\eprp}{\end{prp}}

\newcommand{\bfct}{\begin{fct}}
\newcommand{\efct}{\end{fct}}
\newtheorem{prob}[thm]{Problem}
\newcommand{\bprb}{\begin{prob}}
\newcommand{\eprb}{\end{prob}}
\newtheorem{lem}[thm]{Lemma}
\newcommand{\blem}{\begin{lem}}
\newcommand{\elem}{\end{lem}}
\newtheorem{cor}[thm]{Corollary}
\newcommand{\bcor}{\begin{cor}}
\newcommand{\ecor}{\end{cor}}
\newtheorem{conj}[thm]{Conjecture}
\newcommand{\bcnj}{\begin{conj}}
\newcommand{\ecnj}{\end{conj}}
\theoremstyle{definition}
\newtheorem{defn}[thm]{Definition}
\newcommand{\bdfn}{\begin{defn}}
\newcommand{\edfn}{\end{defn}}
\theoremstyle{remark}
\newtheorem{rem}[thm]{Remark}
\newcommand{\brem}{\begin{rem}}
\newcommand{\erem}{\end{rem}}

\newcommand{\bexm}{\begin{exm}}
\newcommand{\eexm}{\end{exm}}
\newcommand{\bpf}{\begin{proof}}
\newcommand{\epf}{\end{proof}}

\newcommand{\be}{\begin{enumerate}}
\newcommand{\ee}{\end{enumerate}}
\newcommand{\bi}{\begin{itemize}}
\newcommand{\itm}{\item}
\newcommand{\ei}{\end{itemize}}
\newcommand{\bdesc}{\begin{description}}
\newcommand{\edesc}{\end{description}}
\forget
\forgotten

\newcommand{\sone}{\mathsf{S}_1}
\newcommand{\sfin}{\mathsf{S}_{fin}}
\newcommand{\ufin}{\mathsf{U}_{fin}}

\author{Boaz Tsaban}
\thanks{Supported by the Koshland Center for Basic Research.}
\address[Boaz Tsaban]{Department of Mathematics,
Weizmann Institute of Science,
Rehovot 76100,
Israel; and
Department of Mathematics,
Bar-Ilan University,
Ramat-Gan 52900,
Israel.}
\email{boaz.tsaban@weizmann.ac.il}
\urladdr{http://www.cs.biu.ac.il/\~{}tsaban}

\author{Lyubomyr Zdomskyy}
\address[Lyubomyr Zdomskyy]{Department of Mechanics and Mathematics,
Iv\-an Franko Lviv National University, Universytetska 1, Lviv
79000, Ukraine; and Department of Mathematics, Weizmann
Institute of Science, Rehovot 76100, Israel.}
\email{lzdomsky@rambler.ru}



\title[On the Pytkeev property]{On the Pytkeev property in spaces of continuous functions (II)}

\begin{document}

\begin{abstract}
We prove that for each Polish space $X$, the space $C(X)$ of continuous
real-valued functions on $X$ satisfies (a strong version of) the Pytkeev
property, if endowed with the compact-open topology.
We also consider the Pytkeev property in the case where
$C(X)$ is endowed with the topology of pointwise convergence.
\end{abstract}

\maketitle

\section{Introduction}

For a topological space $X$,
$C(X)$ is the family of all real-valued continuous
functions on $X$. We consider two standard topologies on $C(X)$,
which make it a topological group.
Let $\zero$ denote the constant zero function on $X$.

$C_k(X)$ denotes $C(X)$, endowed with the compact-open topology.
For a set $K\sbst X$ and $n\in\N$, let
$$\znb{K}{n}=\left\{f\in C_k(X) : (\forall x\in K)\ |f(x)|<\frac{1}{n}\right\}.$$
When $K$ ranges over the compact subsets of $X$ and $n$ ranges over
$\N$, the sets $\znb{K}{n}$ form a local base at $\zero$.

$C_p(X)$ denotes $C(X)$, endowed with the topology of pointwise convergence.
Here, a local base at $\zero$ is given by the sets $\znb{F}{n}$,
where $n\in\N$, and $F$ ranges over the finite subsets of $X$.

$C_k(X)$ is metrizable if, and only if, $X$ is hemicompact
(i.e., there is a countable family of compact sets such that each compact subset of $X$ is
contained in some member of the family) \cite{McCoy80}.
In particular, $C_k(\NN)$ is not metrizable.
Restricting attention to first countable spaces $X$,
McCoy \cite{McCoy80}
observed that for $C_k(X)$ to be metrizable, it suffices that it
has the \emph{Fr\'echet-Urysohn} property, that is, for each
$A\sbst C_k(X)$ with $\zero\in\cl{A}$, there is a sequence of
elements of $A$ converging to $\zero$. Despite the fact that
$C_k(\NN)$ does not have the Fr\'echet-Urysohn property, we show
in Section \ref{CkSec} that it has the slightly weaker Pytkeev
property.

\medskip

As for $C_p(X)$, it is metrizable if, and only if, $X$ is countable \cite{Arhan92}.
Here, the Fr\'echet-Urysohn property does not imply metrizability, and
Sakai asked whether for $C_p(X)$, the Pytkeev property implies the
Fr\'echet-Urysohn property.
We establish several weaker assertions (Section \ref{CpSec}).

\section{The compact-open topology}\label{CkSec}

Let $X$ be a topological space.
$C_k(X)$ has the \emph{Pytkeev property} \cite{Pyt84} if for each $A\sbst C_k(X)$ with
$\zero\in\cl{A}\sm A$, there are infinite sets $A_1,A_2,\dots\sbst A$
such that each neighborhood of $\zero$ contains some $A_n$.

The notion of a \emph{$k$-cover} is central in the study of local properties of $C_k(X)$
(see \cite{AppKII} and references therein).
A cover $\cU$ of $X$ is a $k$-cover of $X$ if $X\nin\cU$, but for each compact $K\sbst X$,
there is $U\in\cU$ such that $K\sbst U$.

\bthm\label{ckbaire}
$C_k(\NN)$ has the Pytkeev property.
\ethm
\bpf
By a theorem of Pavlovic and Pansera \cite{PavPan06},
it suffices to prove that for each open $k$-cover $\cU$ of $X$,
there are infinite sets $\cU_1,\cU_2,\dots \sbst \cU$ such that
$\{\bigcap \cU_n:n\in\N\}$ is a $k$-cover of $X$.
We will show that $\NN$ has the mentioned covering property.

To this end, we set up some basic notation.
For $s\in\N^{<\alephes}$, $[s]=\{f\in\NN : s\sbst f\}$, and $|s|$ denotes
the length of $s$.
For $S\sbst\N^{<\alephes}$, $[S]=\Union_{s\in S}[s]$.
For an open $U\sbst\NN$, $U(n)=\{s\in\N^n : [s]\sbst U\}$.
Note that for each $n$, $[U(n)]\sbst [U(n+1)]$, and $U=\Union_n[U(n)]$.

\blem\label{infbasic}
Assume that $\cU$ is an open $k$-cover of $\NN$. Then:
\be
\itm $\cV=\{[U(n)] : U\in\cU, n\in\N\}$ is a $k$-cover of $\NN$.
\itm There is $n$ such that $\{U(n) : U\in\cU\}$ is infinite.
\itm For each compact $K\sbst\NN$, there is $n$ such that $\{U(n) : U\in\cU, K\sbst [U(n)]\}$ is infinite.
\ee
\elem
\bpf
(1) For each compact $K\sbst\NN$, there is $U\in\cU$ such that $K\sbst U$.
As $U=\Union_n[U(n)]$ and $K$ is compact, there is $n$ such that $K\sbst [U(n)]\in\cV$.

(2) Assume that for each $n$, $\{U(n) : U\in\cU\}$ is finite.
Note that for each $U\in\cU$ and each $n$, $[U(n)]\sbst U\neq\NN$,
and therefore $U(n)\neq\N^n$.
Proceed by induction on $n$:

\smallskip\noindent\emph{Step 1.}
As $\cU(1)=\{U(1) : U\in\cU\}$ is finite and $\N\nin\cU(1)$, there is a finite
$F_1\sbst \N$ which is not contained in any member of $\cU(1)$.

\smallskip\noindent\emph{Step $n$.}
As $\cU(n)=\{U(n) : U\in\cU\}$ is finite and $F_{n-1}\x\N$ is not contained
in any member of $\cU(n)$, there is a finite $F_n\sbst F_{n-1}\x\N$
which is not contained in any member of $\cU(n)$, and such that $F_n\rest(n-1)=F_{n-1}$.

\smallskip
Take $K=\bigcap_n[F_n]$ (the set of all infinite branches through the finitely
splitting tree $\Union_n F_n$). As $K$ is compact, there is $U\in\cU$ such that
$K\sbst U$. As $U=\Union_n[U(n)]$ and $K$ is compact, there is $n$ such that
$K\sbst [U(n)]$. But then $F_n\sbst U(n)$, a contradiction.

\smallskip
(3) By (1), $\{[U(n)] : U\in\cU, n\in\N, K\sbst[U(n)]\}$ is a $k$-cover of $\NN$.
By (2), there is $m$ such that
$$\cV=\{[[U(n)](m)] : U\in\cU, n\in\N, K\sbst[U(n)]\}$$
is infinite.
For all $U$ and $n$, $[[U(n)](m)]$ is equal to $[U(n)]$ when $n\le m$, and to $[U(m)]$ when $m<n$.
Thus, $\cV=\Union_{n\le m}\{[U(n)] : U\in\cU, K\sbst[U(n)]\}$, and therefore there is $n\le m$
such that $\{[U(n)] : U\in\cU, K\sbst[U(n)]\}$ is infinite.
\epf

\newcommand{\LES}[1]{[\le\! #1]}
For each $n$ and $s\in\N^n$, let $\LES{s}=[\{t\in\N^n : t\le s\}]$,
where $\le$ is pointwise.
The following lemma gives more than what is needed
in our theorem.

\blem\label{strongcov}
Let $\cU$ be an open $k$-cover of $\NN$.
There is $S\sbst\N^{<\alephes}$ such that
for each $s\in S$, $\cU_s = \{U\in\cU : \LES{s}\sbst U\}$
is infinite, and $\{\LES{s} : s\in S\}$ is a clopen $k$-cover
of $\NN$ (refining $\{\bigcap\cU_s : s\in S\}$).
\elem
\bpf
We actually prove the stronger result,
that the statement in the lemma holds when
$$\cU_s = \{[U(|s|)] : U\in\cU, \LES{s}\sbst U\}$$
for each $s\in S$.

Let $S$ be the set of all $s\in\N^{<\alephes}$ such that $\cU_s$
is infinite.
If $K\sbst\NN$ is compact, take $f\in\NN$ such that
the compact set $K(f)=\{g\in\NN : g\le f\}$ contains $K$.
By Lemma \ref{infbasic}, there is $n$ such that there are infinitely many
sets $U(n)$, $U\in\cU$, with $K(f)\sbst [U(n)]$, that is, $\LES{f\rest n}\sbst U$.
Thus, $f\rest n\in S$. Clearly, $K\sbst K(f)\sbst \LES{f\rest n}$.
\epf
This completes the proof of Theorem \ref{ckbaire}.
\epf

\bdfn
For shortness, we say that a topological space $X$ is \emph{nice}
if there is a countable family $\cC$ of open subsets of $X$, such that
for each open $k$-cover $\cU$ of $X$,
$\cS=\{V\in\cC : (\exists^\oo U\in\cU)\ V\sbst U\}$ is a $k$-cover
of $X$.
\edfn

By Lemma \ref{strongcov}, $\NN$ is nice.

\bdfn
A topological space $Y$ has the \emph{strong Pytkeev property}
if for each $y\in Y$, there is a \emph{countable} family $\cN$ of subsets of $Y$,
such that for each neighborhood $U$ of $y$ and each
$A\sbst Y$ with $y\in\cl{A}\sm A$, there is $N\in\cN$
such that $N\sbst U$ and $N\cap A$ is infinite.
\edfn

If $Y$ is first countable, then it has the strong Pytkeev property.
The converse fails, even in the realm of $C_k(X)$.
Indeed, $C_k(\NN)$ is not first countable (since it is a non-metrizable
topological group), and we have the following.

\bthm\label{spyt}
$C_k(\NN)$ has the strong Pytkeev property.
\ethm

Theorem \ref{spyt} follows from the following.

\blem
If $X$ is nice, then $C_k(X)$ has the strong Pytkeev property.
\elem
\bpf
Let $\cC$ be as in the definition of niceness for $X$.
It suffices to verify the strong Pytkeev property of $C_k(X)$ at
$\zero$.
Set
$$\cN=\{\znb{V}{n} : V\in\cC, n\in\N\}.$$
Assume that $A\sbst C_k(X)$ and $\zero\in\cl{A}\sm A$.
There are two cases to consider.

Case 1: For each $n$, there is $f_n\in A\cap\znb{X}{n}$ (equivalently, there are infinitely many such $n$).
Given any neighborhood $\znb{K}{m}$ of $\zero$, take $V\in\cC$
with $K\sbst V$. Then $\znb{V}{m}\sbst \znb{K}{m}$,
and $\znb{V}{m}\cap A\spst \{f_n : n\ge m\}$ is infinite.

Case 2: There is $N$ such that for each $n\ge N$, $A\cap\znb{X}{n}=\emptyset$.
Fix $n\ge N$. $\cU_n = \{f\inv[(-1/n,1/n)] : f\in A\}$
is a $k$-cover of $X$.
Thus,
\begin{eqnarray*}
\cS_n & = & \{V\in\cC : (\exists^\oo U\in\cU_n)\ V\sbst U\} \sbst\\
& \sbst & \{V\in\cC : (\exists^\oo f\in A)\ V\sbst f\inv[(-1/n,1/n)]\} =\\
& = & \{V\in\cC : \znb{V}{n}\cap A\mbox{ is infinite}\}
\end{eqnarray*}
is a $k$-cover of $X$.

Consider any (basic) open neighborhood $\znb{K}{n}$ of $\zero$.
Take $V\in\cS_n$ such that $K\sbst V$.
Then $\znb{V}{n}\in\cN$, $\znb{V}{n}\sbst\znb{K}{n}$, and $\znb{V}{n}\cap A$ is infinite.
\epf

A function $f:X\to Y$ is \emph{compact-covering} if for each
compact $K\sbst Y$, there is a compact $C\sbst X$ such that
$K\sbst f[C]$.
Hereditary local properties of a space $C_k(X)$ are clearly preserved when transforming
$X$ by a continuous compact-covering functions.
(Indeed, if $f:X\to Y$ is a continuous compact-covering surjection,
then $g\mapsto g\circ f$ is an embedding of $C_k(Y)$ into $C_k(X)$.)

\bcor
For each Polish space $X$, $C_k(X)$ has the strong Pytkeev property.
\ecor
\bpf
$X$ is the image of $\NN$ under a continuous compact-covering function.
Indeed \cite{Kechris}: There is a closed $C\sbst\NN$ such that
$X$ is the image of $C$ under a perfect (thus compact-covering) function.
As $C$ is closed, it is a retract of $\NN$, and the retraction is clearly
compact covering.
\epf

\section{The topology of pointwise convergence}\label{CpSec}

There is a very rich local-to-global theory, due to \Arh{} and his
followers, which studies local properties of $C_p(X)$ by translating them into covering properties.
An elegant and uniform treatment of covering properties was given by Scheepers \cite{coc1,coc2}.
We recall a part of this theory that puts the results of the present section in their proper
context.

Let $X$ be a topological space. $\cU$ is a \emph{cover} of $X$ if
$X=\Union\cU$ but $X\nin\cU$.
A cover $\cU$ of $X$ is an \emph{$\omega$-cover} of $X$ if for
each finite subset $F$ of $X$, there is $U\in\cU$ such that $F\subseteq U$.
$\cU$ is a \emph{$\gamma$-cover} of $X$ if it is infinite and for each $x$ in
$X$, $x\in U$ for all but finitely many $U\in\cU$.
Let $\cO$, $\Omega$, and $\Gamma$ denote the collections of all open
covers, $\omega$-covers, and $\gamma$-covers of $X$, respectively.
Let $\scrA$ and $\scrB$ be collections of covers of a space $X$.
Following are selection hypotheses which $X$ may satisfy or not
satisfy \cite{coc1}.
\bdesc
\item[$\sone(\scrA,\scrB)$]{
For all $\cU_1,\cU_2,\dots\in\scrA$,
there are $U_1\in\cU_1,U_2\in\cU_2,\dots$, such that $\{U_1,U_2,\dots\}\in\scrB$.
}
\item[$\sfin(\scrA,\scrB)$]{
For all $\cU_1,\cU_2,\dots\in\scrA$, there are finite
$\cF_1\sbst\cU_1,\cF_2\sbst\cU_2,\dots$, such that
$\Union_{n\in\N}\cF_n\in\scrB$.
}
\item[$\ufin(\scrA,\scrB)$]{
For all $\cU_1,\cU_2,\dots\in\scrA$, there are finite
$\cF_1\sbst\cU_1,\cF_2\sbst\cU_2,\dots$, such that
$\{\Union\cF_1,\Union\cF_2,\dots\}\in\scrB$.
}
\edesc

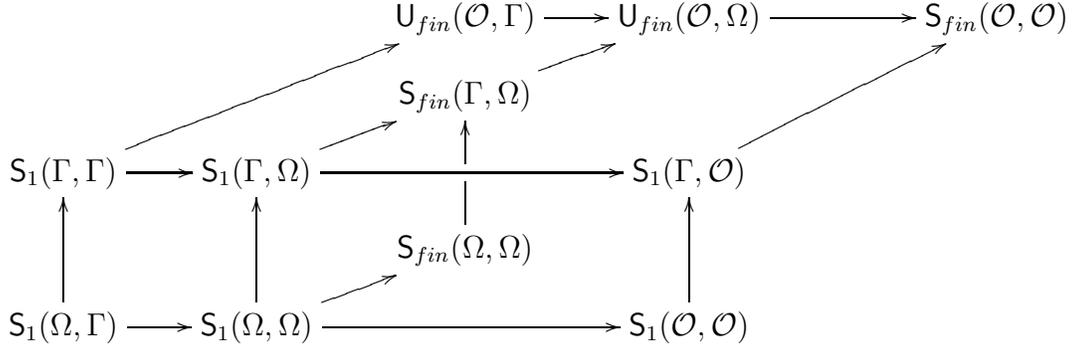
\begin{figure}[!ht]
{%\scriptsize
\begin{changemargin}{-2cm}{-1cm}
\begin{center}
$\xymatrix@R=10pt{%@C=-2pt@R=10pt{%@=7pt{
%1
&
&
& \ufin(\cO,\Gamma)\ar[r]
& \ufin(\cO,\Omega)\ar[rr]
& & \sfin(\cO,\cO)
\\
%2
&
&
& \sfin(\Gamma,\Omega)\ar[ur]
\\
%3
& \sone(\Gamma,\Gamma)\ar[r]\ar[uurr]
& \sone(\Gamma,\Omega)\ar[rr]\ar[ur]
& & \sone(\Gamma,\cO)\ar[uurr]
\\
%4
&
&
& \sfin(\Omega,\Omega)\ar'[u][uu]
\\
%5
& \sone(\Omega,\Gamma)\ar[r]\ar[uu]
& \sone(\Omega,\Omega)\ar[uu]\ar[rr]\ar[ur]
& & \sone(\cO,\cO)\ar[uu]
}$
\caption{The Scheepers Diagram}\label{SchDiagram}
\end{center}
\end{changemargin}
}
\end{figure}

Some of the properties defined in this manner
were studied earlier by Hurewicz ($\ufin(\cO,\Gamma)$), Menger ($\sfin(\cO, \cO)$),
Rothberger ($\sone(\cO, \cO)$, traditionally known as the $C''$ property),
Gerlits and Nagy ($\sone(\Omega,\Gamma)$, traditionally known as the $\gamma$-property),
and others.
Each of these properties is either trivial, or equivalent to
one in Figure \ref{SchDiagram} (where an arrow denotes implication)
\cite{coc2}.

In the remainder of this paper, all spaces $X$ are assumed to be Tychonoff.
A space $X$ satisfies $\sone(\Omega,\Gamma)$ if, and only if,
$C_p(X)$ has the Fr\'echet-Urysohn property \cite{GN}. In particular,
if $X$ satisfies $\sone(\Omega,\Gamma)$, then $C_p(X)$ has the Pytkeev property.
\bprb[Sakai \cite{Sakai06}]\label{sap}
Assume that $C_p(X)$ has the Pytkeev property. Must $X$ satisfy $\sone(\Omega,\Gamma)$?
\eprb
For metric spaces $X$ which are countable unions of totally bounded subspaces,
Miller proved that consistently, $X$ is countable whenever $C_p(X)$ has the Pytkeev property
(this is essentially proved in Theorem 18 of \cite{Pyt}).
It follows that a positive answer to Sakai's Problem \ref{sap} is consistent in this realm.
However, we suspect that the following holds.

\bcnj[\textsf{CH}]
There is $X\sbst\NN$ such that $C_p(X)$ has the Pytkeev property,
but $X$ does not even satisfy Menger's property $\sfin(\cO,\cO)$.
\ecnj

It is therefore natural to consider the conjunction of ``$C_p(X)$ has the Pytkeev property''
with properties in the Scheepers Diagram \ref{SchDiagram}.

A combination of results of Ko\v{c}inac and Scheepers \cite{coc7} and Sakai \cite{Sakai06}
gives that if
$C_p(X)$ has the Pytkeev property and $X$ satisfies $\sfin(\Omega,\Omega)$, then all finite
powers of $X$ satisfy $\ufin(\cO,\Gamma)$ as well as $\sone(\cO,\cO)$.
We will prove several results of a similar flavor.

The combinatorial terminology in the remainder of the paper is as follows:
For $f,g\in\NN$, $f\le^* g$ means $f(n)\le g(n)$ for all but finitely many $n$.
$B\sbst\NN$ is \emph{bounded} if there is $g\in\NN$ such that for each $f\in B$,
$f\le^* g$.
$D\sbst\NN$ is \emph{finitely dominating} if its closure under pointwise maxima of
finite subsets is dominating.

\bthm\label{OOm}
If $C_p(X)$ has the Pytkeev property and $X$ satisfies $\ufin\allowbreak(\cO,\Omega)$, then
$X$ satisfies $\ufin(\cO,\Gamma)$ as well as $\sone(\cO,\cO)$.
\ethm
\bpf
As $C_p(X)$ has the Pytkeev property,
$X$ is Lindel\"of and zero-dimensional \cite{Sakai03}.
This is needed for the application of the quoted combinatorial
theorems below.

We first prove that $X$ satisfies $\ufin(\cO,\Gamma)$.
By \cite{Rec94}, it suffices to prove the following.
\blem\label{bounded}
If $C_p(X)$ has the Pytkeev property and $X$ satisfies $\ufin\allowbreak(\cO,\Omega)$, then
each continuous image $Y$ of $X$ in $\NN$ is bounded.
\elem
\bpf
Let $Y$ be a continuous image of $X$ in $\NN$.
Since we can transform $Y$ continuously by $f(n)\mapsto f(0)+f(1)+\dots+f(n)+n$,
we may assume that all elements of $Y$ are increasing.
If there is an infinite $I\sbst\N$ such that $\{f\rest I : f\in Y\}$ is
bounded, then $Y$ is bounded. We therefore assume that there is $N$
such that for each $n\ge N$, $\{f(n) : f\in Y\}$ is infinite.

As $Y$ satisfies $\ufin(\cO,\Omega)$, $Y$ is not finitely dominating \cite{huremen1},
that is, there is $g\in\NN$ such that the clopen sets $U_n = \{f\in Y : f(n)\le g(n)\}$,
$n\ge N$, form an $\omega$-cover of $Y$. As $C_p(Y)$ has the Pytkeev property,
there are infinite $I_1,I_2,\dots\sbst\N\sm\{0,\ldots,N-1\}$
such that $\{\bigcap_{k\in I_n} U_k : n\in\N\}$ is an $\w$-cover of $Y$ \cite{Sakai03}.
For each $n$, $\{f\rest I_n : f\in \bigcap_{k\in I_n} U_k\}$ is bounded,
and therefore $\bigcap_{k\in I_n} U_k$ is bounded. Thus, $Y=\Union_n\bigcap_{k\in I_n} U_k$ is bounded.
\epf

We now show that $X$ satisfies $\sone(\cO,\cO)$.
It suffices to prove that each continuous image $Y$ of $X$ in $\NN$ has
strong measure zero with respect to the standard metric of $\NN$ \cite{FM}.
Indeed, by Lemma \ref{bounded}, such an image $Y$ is bounded, and thus is a countable union of totally
bounded subspaces of $\NN$. By a theorem of Miller \cite{Pyt}, if $C_p(Y)$ has the Pytkeev property
and $Y$ is a countable union of totally bounded subspaces, then $Y$ has strong measure zero.
\epf

$\Dfin$ is the family of all subsets of $\NN$ which are not
finitely dominating, and $\cov(\Dfin)=\min\{|\cF| : \cF\sbst\Dfin\mbox{ and }\bigcup\cF=\NN\}$.
The hypothesis $\cov(\Dfin)<\fd$ holds, e.g., in the Cohen reals model,
or if $\fd$ is singular \cite{ShTb768}.

\bthm[$\cov(\Dfin)<\fd$]
Assume that for each $Y\sbst X$, $C_p(Y)$ has the Pytkeev property.
Then $X$ satisfies $\ufin(\cO,\Gamma)$ as well as $\sone(\cO,\cO)$.
\ethm
\bpf
By Theorem \ref{OOm}, it suffices to prove that $X$ satisfies $\ufin(\cO,\Omega)$,
or equivalently, that no continuous image $Y$ of $X$ in $\NN$ is
finitely dominating.

Assume that $Y$ is a continuous image of $X$ in $\NN$.
We may assume that all elements of $Y$ are increasing.
Let $\kappa=\cov(\Dfin)<\fd$, and $Y_\alpha\sbst\NN$, $\alpha<\kappa$, be not finitely
dominating and such that $\Union_{\alpha<\kappa}Y_\alpha = \NN$.
For each $\alpha<\kappa$, $Y\cap Y_\alpha$ is not finitely dominating,
and since it is a continuous image of a subset of $X$, $C_p(Y\cap Y_\alpha)$
has the Pytkeev property.
The proof of Lemma \ref{bounded} shows the following.
\blem
Assume that $Z\sbst\NN$, all elements of $Z$ are increasing,
$Z$ is not finitely dominating, and $C_p(Z)$ has the Pytkeev property.
Then $Z$ is bounded.\hfill\qed
\elem
It follows that $Y\cap Y_\alpha$ is bounded for all $\alpha<\kappa$, and as $\kappa<\fd$,
$Y=\Union_{\alpha<\kappa}Y\cap Y_\alpha$ is not finitely dominating.
\epf

We now consider the \emph{strong} Pytkeev property of $C_p(X)$.
A space $Y$ has a \emph{countable cs$^*$-character} \cite{BZ04} if for each
$y\in Y$, there is a countable family $\cN$ of subsets of $Y$,
such that for each sequence in $Y$ converging to $y$ (but not eventually
equal to $y$) and each neighborhood $U$ of $y$, there is $N\in\cN$ such that
$N\sbst U$ and $N$ contains infinitely many elements of that sequence.
Clearly, the strong Pytkeev property implies countable cs$^*$-character.
For topological groups, the conjunction of countable cs$^*$-character and the
Fr\'echet-Urysohn property implies metrizability \cite{BZ04}.
As $C_p(X)$ is a topological group, we have the following.

\bcor\label{stronger}
If $C_p(X)$ has the Fr\'echet-Urysohn property as well as the strong Pytkeev property,
then $X$ is countable.\hfill\qed
\ecor

As the Pytkeev property follows from the Fr\'echet-Urysohn property, we
have the following.

\bcor
The Pytkeev property for $C_p(X)$ does not imply the strong Pytkeev property
for $C_p(X)$.\hfill\qed
\ecor

If, consistently, there is an uncountable $X$ such that
$C_p(X)$ has the strong Pytkeev property, then the answer to Sakai's Problem \ref{sap}
is negative: By corollary \ref{stronger}, in this case $C_p(X)$
cannot have the Fr\'echet-Urysohn property.\footnote{Unfortunately,
this strategy does not work: Sakai has recently proved
that if $C_p(X)$ has the strong Pytkeev property (or even just countable cs$^*$-character),
then $X$ is countable \cite{Sakaics*}. This extends Corollary \ref{stronger}, and can be contrasted
with Theorem \ref{ckbaire}.}

\subsubsection*{Acknowledgement} We thank Masami Sakai for his useful comments.

\ed
\begin{thebibliography}{00}

\bibitem{Arhan92}
A.\ V.\ Arhangel'ski\u{\i},
\textbf{Topological Function Spaces},
Kluwer Academic Publishers, 1992.

\bibitem{BZ04}
T.\ Banakh and L.\ Zdomskyy,
\emph{The topological structure of (homogeneous) spaces and groups with countable $\mathrm{cs}^\ast$-character},
Applied General Topology \textbf{5} (2004), 25--48.

\bibitem{AppKII}
A. Caserta, G. DiMaio, Lj.\ D.R. Ko\v{c}cinac, E. Meccariello,
\emph{Applications of $k$-covers II},
Topology and its Applications \textbf{153} (2006), 3277--3293.

\bibitem{FM}
D.\ H.\ Fremlin and A.\ W.\ Miller,
\emph{On some properties of Hurewicz, Menger and Rothberger},
Fundamenta Mathematica \textbf{129} (1988), 17--33.

\bibitem{GN}
J.\ Gerlits and Zs.\ Nagy,
\emph{Some properties of $C(X)$, I},
Topology and its Applications \textbf{14} (1982), 151--161.

\bibitem{coc2}
W.\ Just, A.\ W.\ Miller, M.\ Scheepers, and P.\ J.\ Szeptycki,
\emph{The combinatorics of open covers II},
Topology and its Applications \textbf{73} (1996), 241--266.

\bibitem{Kechris}
A.\ S.\ Kechris,
\textbf{Classical Descriptive Set Theory},
Graduate Texts in Mathematics \textbf{156}, Springer-Verlag, 1994.

\bibitem{coc7}
Lj.\ D.R.\  Ko\v{c}inac and M.\ Scheepers,
\emph{Combinatorics of open covers (VII): Groupability},
Fundamenta Mathematicae \textbf{179} (2003), 131--155.

\bibitem{McCoy80}
R. A. McCoy,
\emph{Function spaces which are $k$-spaces},
Topology Proceedings \textbf{5} (1980), 139--146.

\bibitem{PavPan06}
B. Pansera, V. Pavlovi\'c,
\emph{Open covers and function spaces},
Matematicki Vesnik \textbf{58} (2006), 57--70.

\bibitem{Pyt84}
E.\ G.\ Pytkeev,
\emph{On maximally resolvable spaces},
Proceedings of the Steklov Institute of Mathematics \textbf{154} (1984), 225--230.

\bibitem{Rec94}
I.\ Rec\l{}aw,
\emph{Every Luzin set is undetermined in the point-open game},
Fundamenta Mathematicae \textbf{144} (1994), 43--54.

\bibitem{Sakai03}
M.\ Sakai,
\emph{The Pytkeev property and the Reznichenko property in function spaces},
Note di Matematica \textbf{22} (2003), 43--52.

\bibitem{Sakai06}
M.\ Sakai,
\emph{Special subsets of reals characterizing local properties of function spaces},
in: \textbf{Selection Principles and Covering Properties in Topology} (L.\ D.R.\ Ko\v{c}inac, ed.),
Quaderni di Matematica \textbf{18}, Seconda Universita di Napoli, Caserta.

\bibitem{Sakaics*}
M.\ Sakai,
\emph{Function spaces with a countable cs$^*$-network at a point},
Topology and its Applications \textbf{156} (2008), 117--123.

\bibitem{coc1}
M.\ Scheepers,
\emph{Combinatorics of open covers I: Ramsey theory},
Topology and its Applications \textbf{69} (1996), 31--62.

\bibitem{ShTb768}
S.\ Shelah and B.\ Tsaban,
\emph{Critical cardinalities and additivity properties of combinatorial notions of smallness},
Journal of Applied Analysis \textbf{9} (2003), 149--162.

\Pa{Pyt}{P. Simon and B. Tsaban}{On the Pytkeev property in spaces of continuous functions}{Proceedings of the American Mathematical Society}{136}{2008}{1125}{1135}{\arx{math}{GN}{0606270}}

\bibitem{huremen1}
B.\ Tsaban,
\emph{A diagonalization property between Hurewicz and Menger},
Real Analysis Exchange \textbf{27} (2001/2002), 757--763.

\end{thebibliography}
